\documentclass[A4paper]{article}

\newcommand{\footremember}[2]{%
    \footnote{#2}
    \newcounter{#1}
    \setcounter{#1}{\value{footnote}}%
}
\newcommand{\footrecall}[1]{%
    \footnotemark[\value{#1}]%
} 

\usepackage{style}
\begin{document}

\title{The lowest-order Neural Approximated Virtual Element Method}

\author{S. Berrone\footremember{trailer}{Department of Mathematical Sciences, Politecnico di Torino, Italy
  (stefano.berrone@polito.it, davide.oberto@polito.it, moreno.pintore@polito.it, gioana.teora@polito.it).}, D. Oberto\footrecall{trailer}{}, M. Pintore\footrecall{trailer}{},  G. Teora\footrecall{trailer}{}}
\maketitle

\begin{abstract}
We introduce the Neural Approximated Virtual Element Method, a novel polygonal method that relies on neural networks to eliminate the need for projection and stabilization operators in the Virtual Element Method. In this paper, we discuss its formulation and detail the strategy for training the underlying neural network. The efficacy of this new method is tested through numerical experiments on elliptic problems.
\end{abstract}


\section{Introduction}

In the last years, the Virtual Element Method (VEM), introduced in \cite{basicVEM}, has gained considerable interest for its flexibility in handling arbitrary geometries and for the simplicity to build high-order methods. This versatility stems from the introduction of some non-polynomial functions in the local space and the careful choice of the local degrees of freedom that allow to exactly compute suitable polynomial projections of these functions, while maintaining their real nature pure virtual. However, these virtual functions also represent the main limitation of the VEM: the requirement to compute the polynomial projection of the virtual functions to access their point-wise evaluation leads not only to many limitations in a post-processing phase \cite{Credali2023reduced}, but also to a lack of coercivity of the method. Indeed, their introduction entails the need to define a proper stabilization term to enforce the coercivity of the local discrete bilinear form, which is highly problem-dependent \cite{russo2023quantitative} and whose isotropic nature can limit the accuracy of the method in case of strongly anisotropic problems \cite{Gioana_Francesca, borio2023lowest}.  

Very recently, some methods have been proposed in \cite{Credali2023reduced, trezzi2023rational} to approximate these functions and address these drawbacks. However, they further increase the additional computational cost already required by VEM with respect to standard procedures like the Finite Element Method (FEM). 

Inspired by the recent success of Scientific Machine Learning  \cite{cuomo2022scientific}, in this paper we introduce the \textit{Neural Approximated Virtual Element Method} (NAVEM), a novel method in which neural networks are used to approximate the unknown VEM basis functions. Through such approximations, the need for computing the local projection matrices and defining a stability operator is circumvented, thereby aligning NAVEM with a FEM on polygonal meshes.

The paper is organized as follows. In Section \ref{sec:vem}, we briefly introduce the VEM formulation, which is crucial to devise a proper architecture and training strategy for the neural network. The new method is presented in Section \ref{sec:elemental_model} and tested in Section \ref{sec:numerical_results}.

\section{The Virtual Element Method}\label{sec:vem}

Let $\Omega \subset \R^2$ be an open convex bounded polygonal domain and let $\Gamma = \partial \Omega$ be its boundary. Let $\Th$ be a decomposition of $\Omega$ into polygons $E$ and let $\Eh[E]$ be the set of edges of the element $E\in\Th$. Furthermore, we denote by $\Nv[E]$ the number of vertices (and of edges) of the element $E$. Here, as usual, $h$ denotes the maximum diameter of the polygons in $\Th$. We assume that the following mesh assumptions hold true.

\begin{assum}[Mesh assumptions]\label{ass:mesh}
There exists a positive constant $\rho$, independent of
$E$ and $h$, such that
\begin{itemize}
\item each polygon $E \in \Th$ is star-shaped with respect to a ball of radius $\geq \rho h_E$;
\item for each edge $e \in \Eh[E]$, it holds: $\vert e \vert \geq \rho h_E$.
\end{itemize}
\end{assum}

Given a polygon $E$, for each integer $k \geq 0$, we denote by $\Poly{k}{E}$ the set of two-dimensional polynomials of degree up to $k$ defined on $E$, with dimension $n_k = \dim \Poly{k}{E} = \frac{(k+1)(k+2)}{2}$. Furthermore, we introduce 
the set 
\begin{equation*}
\Bk{1}{\partial E} = \left\{v \in C^0\left(\partial E\right) : v_{|e}\in \Poly{1}{e} \forall e \in \Eh[E]\right\},
\end{equation*}
whose dimension is $\dim \Bk{1}{\partial E} = \Nv[E]$.
For all $E \in \Th$, we introduce the lowest-order local virtual element space \cite{basicVEM}
\begin{align}
\label{eq:lapl} \Vh[E]{1} = \Big\{ v \in \sob{1}{E}:\ &(i)\ \Delta v = 0, \\ \label{eq:boundary}
&(ii) \ v_{| \partial E} \in \Bk{1}{\partial E}\Big\}, 
\end{align}
with dimension $\Ndof[E] = \dim \Vh[E]{1} = \Nv[E]$, and we consider the value of $v_h \in \Vh[E]{1}$ at the vertices as degrees of freedom.

\begin{remark}\label{rem:pol}
We observe that $\Poly{1}{E} \subseteq \Vh[E]{1}$, and, if $N_E^v>3$, then $\Poly{1}{E} \subsetneq \Vh[E]{1}$.
\end{remark}

\subsection{The model problem and its virtual element discretization}

Given $f \in \leb{2}{E}$, we consider the Poisson problem with homogeneous Dirichlet boundary conditions: 
\begin{equation}
\begin{cases}
- \Delta  u  = f & \text{in } \Omega, \\
u = 0 & \text{on } \Gamma.
\end{cases}
\label{eq:primalcontinuous}
\end{equation} 
The variational formulation of problem \eqref{eq:primalcontinuous} reads: \textit{Find $u\in \V = \sob[0]{1}{\Omega}$ such that}
\begin{equation}
\dbilin{u}{v} = \sum_{E \in \Th} \dbilin[E]{u}{v} = \scal[\Omega]{f}{v}, \quad \forall v \in \V, 
\label{eq:varcontinuous}
\end{equation} 
where the local bilinear form $\dbilin[E]{}{}$ is defined as
\begin{equation}
\dbilin[E]{u}{v} = \int_E \nabla u \cdot \nabla v, \quad \forall u,v \in \V. 
\label{eq:bilin_cont}
\end{equation}
Now, we can define the lowest-order global virtual element space as
\begin{equation}
\Vh{1} = \{v \in \sob[0]{1}{\Omega} \cap C^0(\overline{\Omega}): v_{|E} \in \Vh[E]{1} \ \forall E \in \Th \}.
\end{equation}
We introduce the local polynomial projectors
\begin{itemize}
    \item $\proj{\nabla}{1} : \V \to \Poly{1}{\Th}$ such that, for each $E\in \Th$,
    \begin{equation*}
        \scal[E]{\nabla v - \nabla \proj{\nabla}{1} v}{\nabla p} = 0,\quad \forall p \in \Poly{1}{E} \text{ and } \int_{\partial E} \proj{\nabla}{1} v = \int_{\partial E} v;
    \end{equation*}
    \item $\proj{0}{0} : \V \to \Poly{0}{\Th}$ which is locally defined as the $\leb{2}{E}$-projector onto constants.
\end{itemize}
Finally, the virtual element discretization of problem \eqref{eq:varcontinuous} reads: \textit{Find $u_h \in \Vh{1}$ such that} 
\begin{equation}
\sum_{E \in \Th} \dbilinh[E]{u_h}{v_h} = \sum_{E \in \Th} \scal[E]{f}{\proj{0}{0} v_h} \quad \forall v_h \in \Vh{1},
\label{eq:vemvardiscrete}
\end{equation} 
where the discrete bilinear form $\dbilinh[E]{}{}$ is the sum of a consistency term related to the accuracy and of a stabilization term $\stab[E]{}{}$ enforcing the coercivity, that is
\begin{equation}
\dbilinh[E]{u_h}{v_h} = \int_E \nabla \proj{\nabla}{1} u_h \cdot \nabla \proj{\nabla}{1} v_h + \stab[E]{(I-\proj{\nabla}{1})u_h}{(I-\proj{\nabla}{1})v_h}.
\end{equation}
In particular, the stabilization term $\stab[E]{}{}$ can be chosen as any symmetric positive definite bilinear form satisfying the following property \cite{basicVEM}: there exist two positive constants $\alpha_{\ast},\ \alpha^{\ast}$ independent of $h$ such that
\begin{equation*}
    \alpha_{\ast} \dbilin[E]{v}{v} \leq \stab[E]{v}{v} \leq \alpha^{\ast} \dbilin[E]{v}{v},\quad \forall v \in \Vh[E]{1}:\ \proj{\nabla}{1}v = 0.
\end{equation*}

\section{The Neural Approximated Virtual Element Method}\label{sec:elemental_model}

Let us introduce the set of the VEM Lagrangian basis functions $\{\varphi_{i}\}_{i=1}^{\Ndof}$ corresponding to the aforementioned degrees of freedom, each of them associated with a different internal vertex $v_i$ of the tessellation $\Th$. We denote by $\mathbb{S}_i = \supp{\varphi_i} = \cup_{j=1}^{N_{v_i}} E_j$ the support of $\varphi_i$, which coincides with the union of the $N_{v_i}$ elements $E_j \in \Th$ adjacent to the vertex $v_i$. Furthermore, given an element $E\in \Th$, for the sake of brevity, we indicate by $\{\varphi_{j,E}\}_{j=1}^{\Ndof[E]}$ the set of the restrictions to $E$ of the Lagrangian basis functions related to the vertices of $E$. Clearly, the local and the global virtual element spaces can be written as
\begin{equation*}
    \Vh[E]{1} = \myspan \{ \varphi_{j,E}: \ j = 1,\dots,\Ndof[E]\}
\end{equation*}
and
\begin{equation*}
    \Vh{1} = \myspan \{ \varphi_{i}: \ i = 1,\dots,\Ndof\}.
\end{equation*}

\subsection{The approximated basis functions}
In this section we describe how to construct an approximated space $\nVh[]{1}$ of $\Vh[]{1}$. 

Let us introduce the space $\HPoly{\ell}{E}$ of the harmonic polynomials of degree up to $\ell \geq 1$ with dimension $\dim \HPoly{\ell}{E} = 2 \ell + 1$. Now, we introduce a neural network which aims to learn the highly non-linear map that associates
\begin{equation}
    (v_j, E) \mapsto \nvarphi_{j,E} \in \HPoly{\ell}{E}, \text{ for each vertex $v_j$ of $E$ \text{ and } $\forall E \in \Th$,}
    \label{eq:nn_map}
\end{equation}
where the functions $\nvarphi_{j,E}$ aim to locally approximate the Lagrangian VEM basis functions. To achieve this goal, firstly we must require that the functions $\nvarphi_{j,E}$ belong to the VEM space $\Vh[E]{1}$, and, in particular, that they satisfy properties \eqref{eq:lapl} and \eqref{eq:boundary} that characterize this space.

In this regard, we note that property \eqref{eq:lapl} is trivially satisfied by functions $\nvarphi_{j,E}$ by construction. Instead, property \eqref{eq:boundary} is, in general, not satisfied by functions belonging to $\HPoly{\ell}{E}$. Nevertheless, we overcome this issue by training the neural network to learn functions $\nvarphi_{j,E}$ in such a way that they mimic the Lagrangian basis functions $\varphi_{j,E}$ at the boundary of the element $E$, where all the virtual functions are known in a closed form. By proceeding in this way, the NAVEM functions $\nvarphi_{j,E}$ approximate the VEM Lagrangian basis functions, i.e.
\begin{equation*}
    \nvarphi_{j,E} \approx \varphi_{j,E} \text{ on } E, \forall j =1,\dots,\Ndof[E], \text{ and } \forall E \in \Th.
\end{equation*}

Now, for each internal vertex $v_i$ of the tessellation $\Th$, we can piece-wise define the functions
\begin{equation*}
    \nvarphi_i = 
\begin{cases}
    \nvarphi_{j,E} & \text{if $v_i$ is the $j$-the vertex of $E$ and $E \in \mathbb{S}_i$},  \\
    0 & \text{otherwise},
\end{cases}
\end{equation*}
and define the local and global lowest-order NAVEM spaces as the sets
\begin{equation*}
\nVh[E]{1} = \myspan\{\nvarphi_{j,E},\ j=1,\dots,\Ndof[E]\} \subset \HPoly{\ell}{E},
\end{equation*}
\begin{equation*}
    \nVh{1} = \myspan\{\nvarphi_{i},\ i=1,\dots,\Ndof\},
\end{equation*}
which approximates the corresponding VEM spaces by construction.

\begin{remark}
We note that we decide to employ the \textit{enlarged} space $\HPoly{\ell}{E}$ of harmonic polynomials of degree up to $\ell \geq 1$ for the approximation of VEM basis functions since it is not possible to accurately approximate any arbitrary element of $\Vh[E]{1}$ with an element of $\Poly{1}{E}$, as highlighted in Remark \ref{rem:pol}. By proceeding in this way, we observe that, with appropriate training, we can approximate the Virtual Element basis functions over all the triangles and rectangles since we are learning functions belonging to a space that contains $\Poly{1}{E} = \myspan \{1,x,y\}$ and, if we choose $\ell \geq 2$, also $xy$.
\end{remark}

Finally, the NAVEM discretization of problem \eqref{eq:varcontinuous} reads: \textit{Find $u_h^{\NN} \in \nVh{1}$ such that}
\begin{equation}
\sum_{E \in \Th} \dbilin[E]{u_h^{\NN}}{v_h^{\NN}} = \sum_{E \in \Th} \scal[E]{f}{v_h^{\NN}} \quad \forall v_h^{\NN} \in \nVh{1},
\label{eq:nnvarcontinuous}
\end{equation} 
where the bilinear form $\dbilin[E]{}{}$ is defined in \eqref{eq:bilin_cont}.

\subsection{The neural network}\label{sec:neural_network}

In this section, we show how to train the neural network underlying the NAVEM method (offline phase) and the procedure to predict the NAVEM basis functions (online phase).

\subsubsection{The offline phase}
Let us consider the reference square $\tilde{S}$ centred at the origin with diameter $h_{\tilde{S}}$. Now, we introduce an orthonormal polynomial basis $\{\tilde{p}_{\beta}\}_{\beta=1}^{2\ell + 1}$ for $\HPoly{\ell}{\tilde{S}}$, which is built by orthogonalizing the scaled polynomial basis for $\HPoly{\ell}{\tilde{S}}$ defined by the recursive strategy used in \cite{Gioana_Francesca}. More precisely, we orthogonalize this basis by applying the modified Gram-Schmidt algorithm twice to the Vandermonde matrix whose columns contain the evaluations of the scaled polynomials at points defined as nodes of a lattice built over $\tilde{S}$.

\bigskip

Now, for each $\Nv \geq 3$, we perform the following steps:
\begin{enumerate}
    \item We consider a set of randomly generated polygons $\mathcal{S}_{\Nv}$ with $\Nv$ vertices and edges. This process is subject to the restrictions imposed by the mesh assumptions \ref{ass:mesh}.
    \item On each element $E \in \mathcal{S}_{\Nv}$ we perform the mapping $F_E : \srescale{E} \to E$ proposed in \cite{Gioana_Fabio} to reduce the variability of elements seen by the neural network. This mapping is an affine transformation whose definition is based on the inertia tensor of the element $E$. We recall that this mapping generates polygons $\srescale{E}$ with unit diameter and centred in the origin.
    \item For each vertex $v_j$ of $\srescale{E} = F_E^{-1}(E): E \in \mathcal{S}_{\Nv}$, we encode the information representing the pair $(v_j, \srescale{E})$ into a vector $\zz_E \in \R^{N_{0}}$ of a suitable dimension $N_{0}$. More precisely, we consider an affine map $ G_{j,\srescale{E}} : \tilde{E}_j \to \srescale{E}$ that rotates and rescales $\srescale{E}$ to place its $j$-vertex into $(1,0)$. Then, we set $\zz_E = [\tilde{\xx}_1^2,\tilde{\xx}_2^2,\tilde{\xx}_1^3,\tilde{\xx}_2^3,\dots,\tilde{\xx}_1^{\Nv},\tilde{\xx}_2^{\Nv}]$, where $(\tilde{\xx}_1^k,\tilde{\xx}_2^k)$ are the coordinates of the $k$-th vertices of $\tilde{E}_j = G_{j, \srescale{E}}^{-1}(\srescale{E})$ and the vertices of $\tilde{E_j}$ are ordered counter-clock-wise starting from $(\tilde{\xx}_1^1,\tilde{\xx}_2^1) = (1,0)$. Thus, the dimension of the input of the neural network is $N_{0} = 2 (\Nv -1)$. This process is done to reduce the dimension of the input of the neural network and it is possible since the prediction of a basis function is not affected by the prediction of basis functions related to the remaining vertices of the element.
    \item Given an activation function $\sigma$, we define a feed-forward fully-connected neural network with $L$ layers and $N_n$, $n=1,\dots, L$, neurons per layer, i.e. a function $c^\NN:\R^{N_{0}}\rightarrow\R^{2 \ell + 1}$ defined by the recursive expression:
    \begin{equation} 
      \begin{aligned}
     &\zz_n = \sigma({\mathbf A}_n \zz_{n-1} + {\mathbf b}_n), \hspace{1cm} & n = 1,...,L-1, \\
     &\cc^\NN = c^{\NN}(\zz_E) = {\mathbf A}_{L} \zz_{L-1} + b_L &
      \end{aligned}
      \label{eq:nn_main_formula}
      \end{equation}
      where $\zz_0 = \zz_E$ is the input, ${\mathbf A}_n \in\mathbb{R}^{N_n \times N_{n-1}}$ and ${\mathbf b}_n \in\mathbb{R}^{N_n}$, $n=1,...,L$ are the matrices and vectors containing the network weights, while the output $\cc^\NN \in \R^{2 \ell + 1}$ represents the vector of coefficients of functions $\nvarphi_{j,\tilde{E}_j} := \nvarphi_{j,E} \circ F_E \circ G_{j,\srescale{E}}$ with respect to the basis $\{\tilde{p}_{\beta}\}_{\beta=1}^{2\ell + 1}$, i.e. $\nvarphi_{j,\tilde{E}_j}=\sum_{\beta=1}^{2\ell + 1}\cc^\NN_{\beta}(\zz_E) \tilde{p}_{\beta}$.
\end{enumerate}

As mentioned in the previous section, we must train the neural network in such a way that the functions $\nvarphi_{j,E}$ mimic the behaviour of $\varphi_{j,E}$ at the boundary of $E$, where such functions are well-known. Thus, we train each neural network \eqref{eq:nn_main_formula} related to a value $\Nv$ by minimizing the loss function:

\begin{align}
\nonumber \mathcal{L}_{\Nv} = \sum_{E\in \mathcal{S}_{\Nv}} \sum_{j=1}^{\Nv} \sum_{\tilde{\xx}\in X_{\partial \tilde{E}_j}} \Bigg[&
\left( \nvarphi_{j,\tilde{E}_j}(\tilde{\xx
}) - \varphi_{j,\tilde{E}_j}(\tilde{\xx})\right)^2  \\
&+\left( \dfrac{\partial\nvarphi_{j,\tilde{E}_j}}{\partial \tilde{\tt}}(\tilde{\xx}) - \dfrac{\partial\varphi_{j,\tilde{E}_j}}{\partial\tilde{\tt
}}(\tilde{\xx})\right)^4 \Bigg],
\label{eq:loss}
\end{align}
where $X_{\partial \tilde{E}_j}$ denotes the set of control points distributed on $\partial \tilde{E}_j$ and $\dfrac{\partial}{\partial\tilde{\tt}}$ represents the tangential derivative. 

We emphasize that the different exponents used in \eqref{eq:loss} are set in order to obtain two contributions with comparable orders of magnitude.

Finally, it is worth mentioning that the diameter of the reference square $\tilde{S}$ is set in such a way all the elements $\tilde{E}_j$ are contained in $\tilde{S}$. Consequently, the polynomial basis $\{\tilde{p}_{\beta}\}_{\beta=1}^{2\ell + 1}$ is well-scaled for each element $\tilde{E}_j$.

\subsubsection{The online phase}

Now, let us describe the online phase wherein we predict the coefficients of the Lagrangian basis function $\nvarphi_{j,E}$ for each pair $(v_j,E)$ to assemble the linear system related to the NAVEM discretization \eqref{eq:nnvarcontinuous}.

For each element $E\in \Th$, firstly we perform the inertial mapping $\srescale{E} = F_E^{-1}(E)$. Subsequently, for each vertex $v_j$ of $\srescale{E}$:
\begin{enumerate}
\item We encode the information $(v_j,\srescale{E})$ to generate the input $\zz_E$ of the neural network.
\item We predict the coefficients of the corresponding NAVEM basis function $\nvarphi_{j,\tilde{E}_j}$ with respect to the polynomial basis $\{\tilde{p}_{\beta}\}_{\beta=1}^{2\ell + 1}$ using the neural network related to the number of vertices of $E$.
\item We map the function $\nvarphi_{j,\tilde{E}_j}=\sum_{\beta=1}^{2\ell + 1}\cc^\NN_{\beta}(\zz_E) \tilde{p}_{\beta}$ and its gradient  $\tilde{\nabla} \nvarphi_{j,\tilde{E}_j}=\sum_{\beta=1}^{2\ell + 1}\cc^\NN_{\beta}(\zz_E) \tilde{\nabla} \tilde{p}_{\beta}$ back to the original element $E$.
\end{enumerate}
Finally, we can assemble and solve the linear system associated with Problem \eqref{eq:nnvarcontinuous} as in standard FEM solvers, to numerically compute the discrete NAVEM solution $u_h^\NN$.


\section{Numerical results}\label{sec:numerical_results}
In this section, we present some numerical experiments to show the effectiveness of our approach to discretize general second-order elliptic problems. More precisely, we consider the advection-diffusion-reaction problem
\begin{equation}
    \begin{cases}
        \nabla \cdot \left( - \D(\xx) \nabla u\right) + \bbeta(\xx) \cdot \nabla u + \gamma(\xx) u = f & \text{in } \Omega,\\
        u = g_D & \text{on } \Gamma,
    \end{cases}
    \label{eq:diffreacadv}
\end{equation}
where
\begin{equation*}
    \D(\xx) = \begin{bmatrix}
        1 + \xx_2^2 & -\xx_1 \xx_2\\
        -\xx_1 \xx_2 & 1 + \xx_1^2
    \end{bmatrix}, \quad \bbeta(\xx) = \begin{bmatrix}
     \xx_1 \\
     -\xx_2
    \end{bmatrix}, \quad \gamma(\xx) = \xx_1 \xx_2.
\end{equation*}
Furthermore, the forcing term $f$ and the Dirichlet boundary condition $g_D$ are set in such a way the exact solution is
\begin{align}
\nonumber    u(\xx) &= 3 \left((\xx_1 - 0.2) + \frac{\xx_2 - 0.3}{2}\right)^2 + 2 \left( \frac{\xx_1 - 0.7}{2} + (\xx_2 - 0.8) \right)^3 \\
    &\ +\ \sin(2 \pi \xx_1) \sin(3 \pi \xx_2). \label{eq:exact_solution}
\end{align}
The contour line plot of the exact solution is shown in Figure \ref{fig:meshes}.

\subsection{The neural network architecture}
Our neural network is made up of 3 hidden layers of 30 neurons each. We employ the hyperbolic tangent as nonlinear activation functions, initialize the weights using a Glorot normal distribution, and set the polynomial degree $\ell = 5$. 
We use a first-order Adam optimizer for the first 3000 epochs with an exponentially decaying learning rate and successively we switch to a second-order BFGS optimizer until convergence of the loss is reached.

\begin{figure}[h]
\centering
\begin{subfigure}{0.44\textwidth}
    \includegraphics[width=\textwidth] {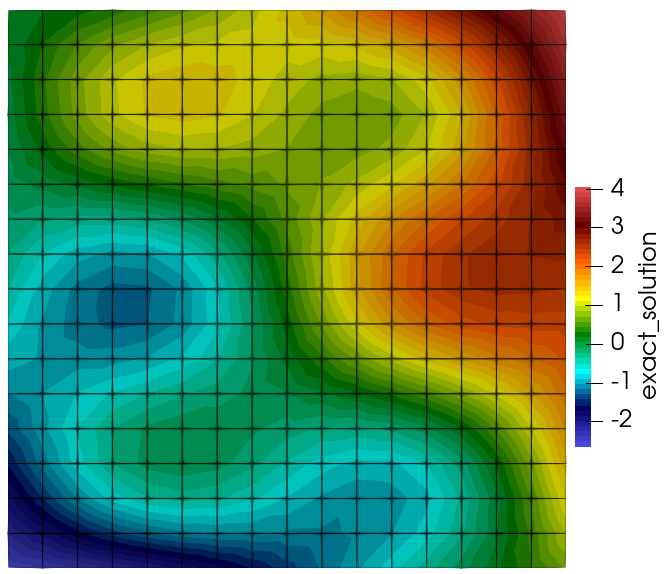}
    \caption{}
\end{subfigure}\quad
\begin{subfigure}{0.44\textwidth}
    \includegraphics[width=\textwidth] {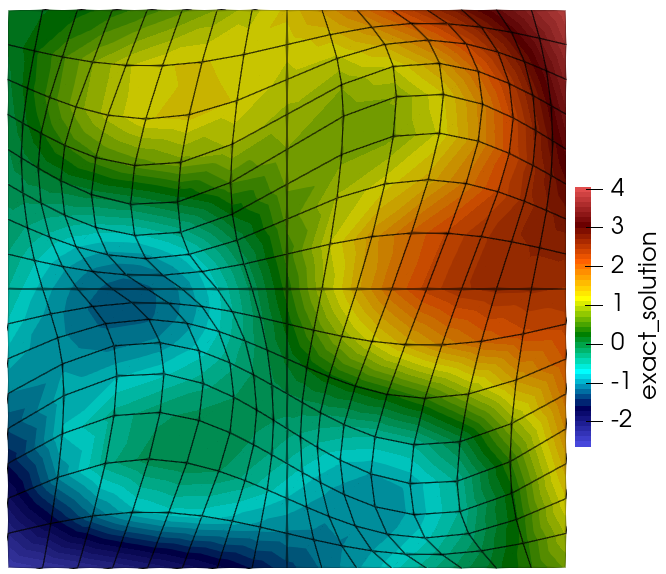}
    \caption{}
\end{subfigure}
\caption{The third mesh of the family coloured by the interpolated exact solution \eqref{eq:exact_solution}. Left: a cartesian mesh. Right: a sine distorted mesh.}
\label{fig:meshes}
\end{figure}

Furthermore, in our experiments, we numerically solve Problem \eqref{eq:diffreacadv} on two different families of meshes made up of quadrilateral elements. The first family consists of $4$ cartesian meshes made up of 16, 64, 256 and 1024 identical elements, respectively. The meshes in the second family are derived from a sine distortion of the related cartesian counterparts. The third refinement of each family is shown in Figure \ref{fig:meshes}.

\begin{figure}[]
\begin{subfigure}{0.3\textwidth}
    \includegraphics[width=\textwidth] {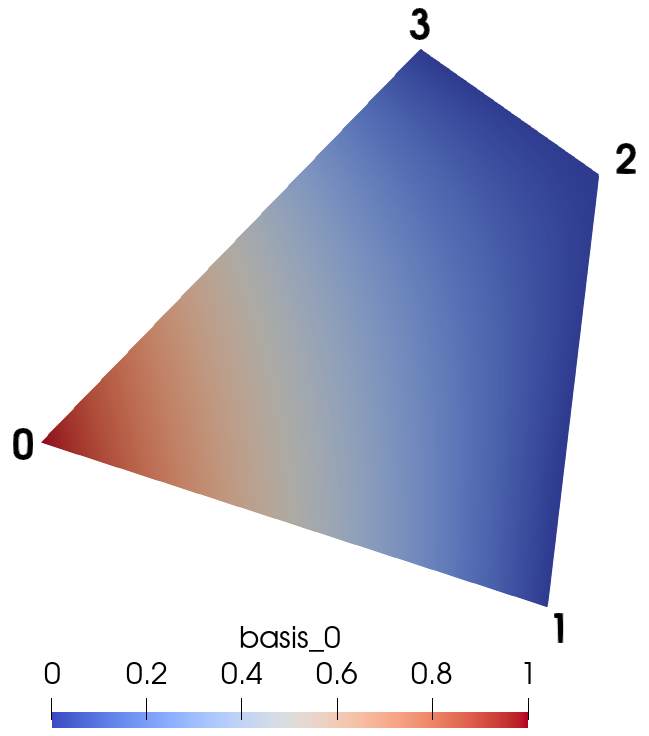}
    \caption{}
\end{subfigure}
\begin{subfigure}{0.3\textwidth}
    \includegraphics[width=\textwidth] {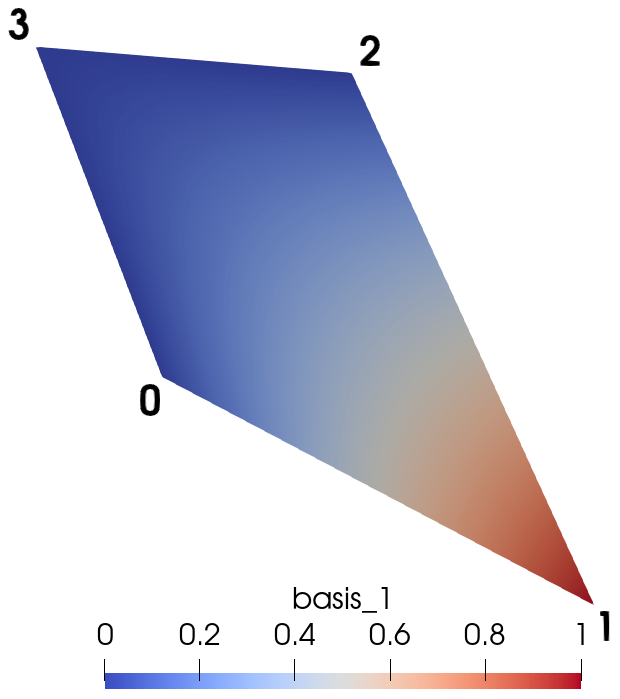}
    \caption{}
\end{subfigure}
\begin{subfigure}{0.3\textwidth}
    \includegraphics[width=\textwidth] {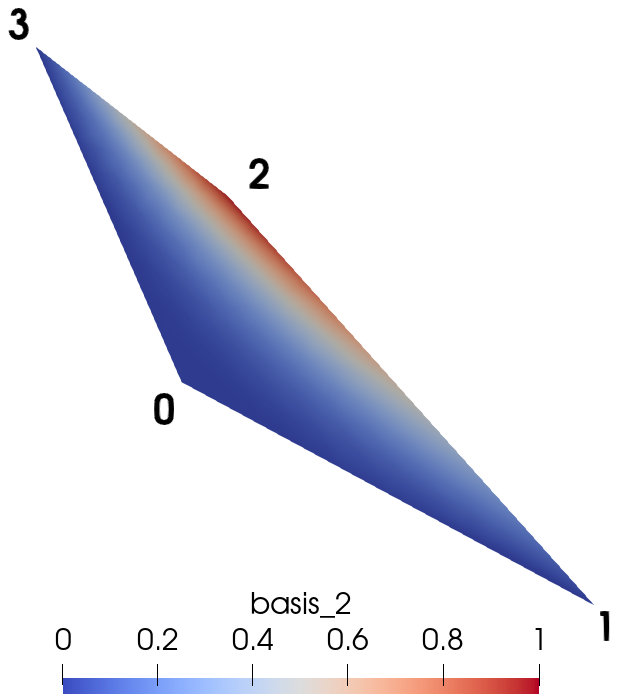}
    \caption{}
\end{subfigure}
\caption{Three different training polygons coloured by the predicted NAVEM basis functions associated with vertices 0, 1 and 2 respectively.}
\label{fig:basis_funcs}
\end{figure}

Since we consider only quadrilateral meshes, we only need to train a single neural network. In this case, the training dataset is made up of 100 randomly generated quadrilaterals. Examples of such generated quadrilaterals coloured by one of the corresponding NAVEM basis functions are shown in Figure \ref{fig:basis_funcs}.

\begin{figure}[]
\begin{subfigure}{0.44\textwidth}
    \includegraphics[width=\textwidth]{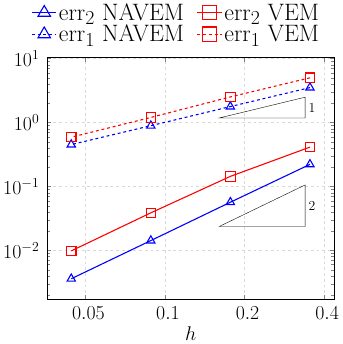}
    \caption{}
\end{subfigure}\quad
\begin{subfigure}{0.44\textwidth}
    \includegraphics[width=\textwidth]{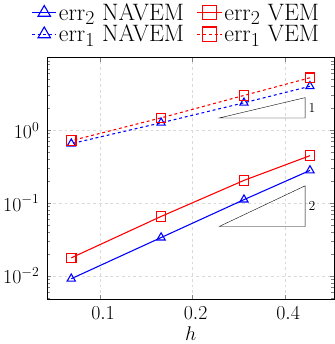}
    \caption{}
\end{subfigure}
\caption{The NAVEM errors \eqref{eq:nn_errors} and the VEM errors \eqref{eq:vem_errors} versus $h$. Left: cartesian meshes. Right: sine distorted meshes.}
\label{fig:errors}
\end{figure}

\subsection{The convergence curves}

To assess the accuracy of our procedure and make a comparison with the standard VEM, we solve Problem \eqref{eq:diffreacadv} with the NAVEM and analyze the behaviour of the following errors
\begin{equation}
    \mathrm{err}_2^2 = \sum_{E\in\Th} \norm[0,E]{u - u^{\NN}_h}^2,\quad \mathrm{err}_1^2 = \sum_{E\in\Th} \norm[0,E]{\nabla u - \nabla u^{\NN}_h}^2,
    \label{eq:nn_errors}
\end{equation}
when $h$ varies.

Since we are dealing with variable coefficients, we opt for a more suitable VEM formulation, detailed in \cite{Beirao2015b}, which is based on the definition of the local enhanced virtual element space and on the introduction of the $L^2$-projector of derivatives. Moreover, since the VEM functions are pure virtual, we are not able to compute the errors as expressed in \eqref{eq:nn_errors}. Thus, in the case of VEM, we define
\begin{equation}
    \mathrm{err}_2^2 = \sum_{E\in\Th} \norm[0,E]{u - \proj{0}{1} u_h}^2,\quad \mathrm{err}_1^2 = \sum_{E\in\Th} \norm[0,E]{\nabla u - \proj{0}{0} \nabla u_h}^2.
    \label{eq:vem_errors}
\end{equation}

The behaviour of the errors \eqref{eq:nn_errors} and \eqref{eq:vem_errors} are shown in Figure \ref{fig:errors}. We observe that the errors \eqref{eq:nn_errors} decay with the same rate of convergence of the related VEM errors, that is $O(h^2)$ for the error in the $L^2$-norm and $O(h^1)$ for the error in the $H^1$-seminorm. Furthermore, we note that the error constants related to NAVEM are smaller than the VEM error constants, that is the curves of convergences of NAVEM are downward shifted with respect to the VEM ones. Since we are approximating the same space, we suppose that a good training of the neural networks could ensure NAVEM achieves also the same theoretical convergence results that hold for the VEM, as supported by the numerical experiments reported here. We highlight that at the moment there exist only a few examples of a priori error estimates for neural networks-based solvers (see for instance \cite{berrone2022variational}). However, in the future, we aim to present a more theoretical discussion.


\section{Conclusions}\label{sec:conclusion}

In this paper, we introduce the Neural Approximated Virtual Element Method. It is a polygonal method that does not require the introduction of any stabilization term and relies on a neural network to approximate the VEM basis functions directly. 

We outline the formulation, describe the neural network training strategy and present some numerical results comparing the performance of NAVEM with the VEM. We observe that we have the same empirical convergence rates achieved by the VEM. Additionally, NAVEM solutions are slightly more accurate than the VEM solutions. This is probably due to the absence of projection or stabilization operators in the NAVEM formulation.

Future perspectives encompass, but are not limited to, the numerical testing of our procedure on more general meshes and the extension of the NAVEM method to higher orders and to the three-dimensional case.

\section*{Acknowledgments}
The author S.B. kindly acknowledges partial financial support provided by PRIN project “Advanced polyhedral discretisations of heterogeneous PDEs for multiphysics problems” (No. 20204LN5N5\_003) and by PNRR M4C2 project of CN00000013 National Centre for HPC, Big Data and Quantum Computing (HPC) (CUP: E13C22000990001). 
The author G.T. kindly acknowledges financial support provided by the MIUR programme ``Programma Operativo Nazionale Ricerca e Innovazione 2014 - 2020'' 
~ (CUP: E11B21006490005). The authors are members of the Italian INdAM-GNCS research group.

\bibliographystyle{IEEEtran}
\bibliography{biblio.bib}

\end{document}